\documentclass[11pt]{article}
\usepackage{fullpage, amsmath, amssymb}

\begin{document}

\title{A Simple Proof of the Quadratic Formula}

\author{Po-Shen Loh\thanks{Department of Mathematical Sciences, Carnegie Mellon
  University. Email: \texttt{po@poshenloh.com}. For a discussion designed for
  the general public, visit \texttt{https://www.poshenloh.com/quadratic}}}

\date{December 16, 2019}

\maketitle

\begin{abstract}
  This article provides a simple proof of the quadratic formula, which also
  produces an efficient and natural method for solving general quadratic
  equations. The derivation is computationally light and conceptually natural,
  and has the potential to demystify quadratic equations for students worldwide.
\end{abstract}

\section{Introduction}

The quadratic formula was a remarkable triumph of early mathematicians, marking
the completion of a long quest to solve quadratic equations, with a storied
history stretching as far back as the Old Babylonian Period around 2000--1600
B.C.  \cite{katz, robson}. 
Over four millennia, many recognized names in mathematics left their mark on
this topic, and the formula became a standard part of a first course in Algebra.

However, it is unfortunate that for billions of people worldwide, the quadratic
formula is also their first experience of a rather complicated formula which
they memorize. Many typically learn it as the systematic alternative to a
guess-and-check method that only factorizes certain contrived quadratic
polynomials. Countless mnemonic techniques abound, from stories of negative bees
considering whether or not to go to a radical party, to songs set to the tune of
\emph{Pop Goes the Weasel.} A derivation by completing the square is usually
included in the curriculum, but its motivation is often challenging for
first-time Algebra learners to follow, and its written execution can be
cumbersome.

This article introduces an independently discovered simple derivation of the
quadratic formula, which also produces a computationally-efficient and natural
method for solving general quadratic equations. The author would actually be
very surprised if this pedagogical approach has eluded human discovery until the
present day, given the 4,000 years of history on this topic, and the billions of
people who have encountered the formula and its proof. Yet this technique is
certainly not widely taught or known. After an earlier version of this arXiv
preprint was publicized by formal and informal media, the author was contacted
by many people with potential references. However, the author still has not
found a previously-existing publicly-shared reference detailing this pedagogical
approach which is mathematically complete and formally correct. (Similar
writings which come close, most notably Savage \cite{savage}, are outlined in
Section \ref{sec:history}.) This article aims to provide a safely referenceable
method and derivation which is logically sound. That said, it is entirely
possible that there is still a previously-existing reference waiting to be
found. So, this article seeks at the very least to popularize a delightful
alternative pedagogical approach to solving quadratic equations, which is
practical for integration into all mainstream curricula.

\section{Derivation}

The reader is encouraged to remember what it was like to be a first-time
Algebra learner, where it took some concentration to combine fractional
expressions, especially those including multiple variables and constants.
This section is intentionally written at that level of simplicity, to
emphasize how straightforward all the algebraic manipulations and concepts
are.

\subsection{Computationally simple derivation of an explicit root formula}
\label{sec:monic}

Throughout this article, we work over the complex numbers. The starting point is
that if we can find a factorization of the following form
\begin{equation}
  \label{eq:monic}
  x^2 + Bx + C = (x-R)(x-S),
\end{equation}
then a value of $x$ makes the product equal zero precisely when at least one of
the factors becomes zero, which happens precisely when $x = R$ or $x = S$. By
the distributive law, it suffices to find two numbers $R$ and $S$ with sum $-B$
and product $C$; then, $\{R, S\}$ will be the complete set of
roots.\footnote{This is the standard method of factoring, which corresponds to
  the converse of relations that are typically attributed to Vi\`ete
  \cite{viete}. Those relations provide inspiration for this method, but are not
  required for logical completeness, because the converse is a straightforward
  consequence of the distributive property. By using the converse, this proof
  does not rely on the theorem that two roots (counting multiplicity) always
exist.}

Two numbers sum to $-B$ precisely when their average is $-\frac{B}{2}$, and so
it suffices to find two numbers of the form $-\frac{B}{2} \pm z$ which multiply
to $C$, where $z$ is a single unknown quantity, because they will automatically
have the desired average.\footnote{This substitution, and this entire solution
  to find two numbers given their sum and product, was known to the Babylonians
  (see, e.g., Burton \cite{burton}, Gandz \cite{gandz-1937}, Irving
  \cite{irving}, or Katz \cite{katz}). It also appeared in the first book of
  Diophantus \cite{diophantus}. This approach therefore represents the fusion of
  these ancient techniques together with Renaissance-era mathematical
  sophistication.  Further historical context follows later in this article.}
(If $z$ turns out to be 0, then we factor with $R = S = -\frac{B}{2}$.) The
product $(-\frac{B}{2} + z)(-\frac{B}{2} - z)$ conveniently matches the form
of a difference of squares, and equals $C$ precisely when
\begin{displaymath}
  \left( -\frac{B}{2} \right)^2 - z^2 = C,
\end{displaymath}
or equivalently, precisely when we have a $z$ which satisfies
\begin{displaymath}
  z^2 = \frac{B^2}{4} - C.
\end{displaymath}
Since a square root always exists (extending to complex numbers if necessary),
arbitrarily select a choice of square root of $\frac{B^2}{4} - C$ to serve as
$z$, in order to satisfy the last equation. Tracing back through the logic, we
conclude that the desired $R$ and $S$ exist in the form $-\frac{B}{2} \pm z$,
and so
\begin{equation}
  \label{eq:formula-monic}
  -\frac{B}{2} \pm \sqrt{\frac{B^2}{4} - C}
\end{equation}
are all the roots of the original quadratic.
\hfill $\Box$

\subsection{Example of use as a method}
\label{sec:example}

The computational and conceptual simplicity of this derivation actually
renders it unnecessary to memorize any formula at all, even for general
coefficients of $x^2$. The proof naturally transforms into a method, and
students can execute its logical steps instead of plugging numbers into a
formula that they do not fully understand. Consider, for example, the
following quadratic:
\begin{displaymath}
  \frac{x^2}{2} - x + 2 = 0.
\end{displaymath}
Multiplying both sides by 2 to make the coefficient of $x^2$ equal to 1, we
obtain the equivalent equation
\begin{displaymath}
  x^2 - 2x + 4 = 0.
\end{displaymath}
If we find two numbers with sum 2 and product 4, then they are all the
solutions. Two numbers have sum 2 precisely when they have average 1.  So, it
suffices to find some $z$ such that two numbers of the form $1 \pm z$ have
product 4 (their average is automatically 1). The final condition is equivalent
to each of these equivalent equations:
\begin{align*}
  1 - z^2 &= 4 \\
  z^2 &= -3.
\end{align*}
We can satisfy the last equation by choosing $i \sqrt{3}$ for $z$. Tracing back
through the logic, we conclude that $1 \pm i \sqrt{3}$ are all the solutions to
the original quadratic. Irrational and imaginary numbers pose no obstacle to
this method.


\subsection{Derivation of traditional quadratic formula with arbitrary
  $\boldsymbol{x^2}$ coefficient}

If one specifically wishes to derive the commonly memorized quadratic
formula using this method, one only needs to divide the equation $ax^2 + bx
+ c = 0$ by $a$ (assume nonzero) to obtain an equivalent equation which matches
the form of \eqref{eq:monic}:
\begin{displaymath}
  x^2 + \left( \frac{b}{a} \right) x + \left( \frac{c}{a} \right) = 0.
\end{displaymath}
Plugging $\frac{b}{a}$ and $\frac{c}{a}$ for $B$ and $C$ in
\eqref{eq:formula-monic}, the roots are:
\begin{displaymath}
  -\frac{b}{2a} \pm \sqrt{\frac{b^2}{4a^2} - \frac{c}{a}}
  \ = \ 
  -\frac{b}{2a} \pm \sqrt{\frac{b^2 - 4ac}{4a^2}}
  \ = \ 
  \frac{-b \pm \sqrt{b^2 - 4ac}}{2a}.
\end{displaymath}
\hfill $\Box$

Observe that with this approach, all of the useful and interesting conceptual
insights are fully isolated in a computationally light derivation of an explicit
formula, while also producing an efficient and understandable algorithm. The
routine but laborious computational portion is required only if a general
formula is sought for memorization purposes.  In light of the efficient
algorithm, however, it becomes questionable whether there is merit to memorizing
a formula without understanding. For example, although the solution to a general
linear equation $ax + b = 0$ is $x = -\frac{b}{a}$ (assume $a \neq 0$), the
equation is typically solved via manipulation instead of plugging into a
memorized formula.

\section{Discussion}

\subsection{Practical relation to other curricular concepts}

Before learning the quadratic formula, students learn how to multiply
binomials, and they see useful expansions such as $(u + v)^2 = u^2 + 2uv +
v^2$ and $(u+v)(u-v) = u^2 - v^2$. Indeed, the first of these expansions is
the cornerstone of the traditional proof of the quadratic formula by
completing the square. The second of these expansions is also of wide
importance: among other things, it is eventually used to rationalize the
denominator of expressions such as $\frac{1}{\sqrt{3} - \sqrt{2}}$ by
multiplying the numerator and denominator by $\sqrt{3} + \sqrt{2}$.

Our approach shows that the factoring method can always be made to work. It
always produces two roots (counting multiplicity) whose sum and product
correspond to coefficients of the quadratic. This therefore presents an
opportunity to prove Vi\`ete's sum and product relations for quadratics.

For first-time Algebra learners, the only new leap of insight is that if one is
seeking two numbers with a desired sum, then they can be parameterized by their
desired average, plus or minus a common unknown amount. In the modern day, a
similar parameterization appears as a useful trick for mentally calculating
products via the difference of squares, such as
\begin{displaymath}
  43 \times 37 = (40 + 3)(40 - 3) = 40^2 - 3^2 = 1591.
\end{displaymath}

This is an ancient trick. Some historians believe the Babylonians used it
thousands of years ago, multiplying in their base-60 number system by
subtracting from tables of squares (see, e.g., Derbyshire
\cite{derbyshire}). It was then natural for them to develop the same
parameterization for finding two numbers, given their sum and product.


\subsection{Comparison to completing the square}

The most common proof of the quadratic formula is via completing the square, and
that was also the method used by al-Khwarizmi \cite{al-khwarizmi} in his
systematic solutions to abstract quadratic equations. Compared to our approach,
the motivation is less direct, as the step of completing the square (for the
simple situation of $x^2 + Bx + C = 0$) simultaneously combines three insights:
\begin{description}
  \item[(i)] The $x^2$ and $Bx$ can be entirely absorbed into a square of
    the form $(x + D)^2$ by using only part of the expansion $(u+v)^2 = u^2
    + 2uv + v^2$ ``backwards,'' to attempt to factor an expression that
    begins with $u^2 + 2uv$.
  \item[(ii)] This perfect square can be created by adding and subtracting the
    appropriate constant, which is $(\frac{B}{2})^2$.
  \item[(iii)] After these manipulations are complete, the equation will have
    $(x + \frac{B}{2})^2$ and some constants, and any such equation can be
    solved by moving constants around and taking a square
    root.\footnote{\label{ftCompSqr} It is interesting to note that this step of
    completing the square uses the fact that the complete set of solutions to
    $x^2 = K$ is $\{\pm \sqrt{K}\}$, which is not obvious to a first-time
  student of quadratics. In contrast, our approach only requires that there
exists an explicit choice of $x$ which satisfies $x^2 = K$, which is often
explicitly constructable.}
\end{description}
The full combination of these insights is required to understand the motivation
for why one should even write down the specific offsetting quantities
$+\frac{B^2}{4} - \frac{B^2}{4}$ in the first line of the completing the square:
\begin{displaymath}
  x^2 + Bx + \frac{B^2}{4} - \frac{B^2}{4} + C = 0.
\end{displaymath}

In contrast, our approach starts from students' existing experience searching
for a pair of numbers with given sum and product, which naturally arises during
the factoring method.  It shows them that the (sometimes frustrating)
guess-and-check process can be replaced by one idea: to parameterize the pair by
its average plus or minus a common unknown offset. No particular formula needs
to be written for the offset itself (unlike the case of carefully selecting
$\frac{B^2}{4}$), and we can simply call it an unknown $z$. Instantly, their
previous experience of trial and error is replaced by a ``forward'' expansion of
the form $(u+v)(u-v) = u^2 - v^2$, which produces an exciting lone $z^2$,
revealing the pair of numbers with all guesswork eliminated.

\subsection{Brief historical context}
\label{sec:history}

Could such a simple proof and pedagogical method possibly be new? The author
researched the English-language literature on the history of mathematics, and
consulted English translations of old manuscripts, from mathematical traditions
ranging from Diophantus \cite{diophantus} to Brahmagupta \cite{brahmagupta},
Yanghui \cite{yanghui}, and al-Khwarizmi \cite{al-khwarizmi}. This section is
too brief to do full justice to the history, and mainly serves to point the
interested reader to relevant resources with much richer detail. In particular,
several books have surveyed the topic of the quadratic formula, such as Chapter
2 of Irving \cite{irving}, and mathematical history books such as Burton
\cite{burton}, Derbyshire \cite{derbyshire}, and Katz \cite{katz}.

As preserved in their cuneiform tablets, the Babylonians had evidence of
formulas for a wide variety of problems of quadratic nature, dating back to the
Old Babylonian Period around 2000--1600 B.C. Although today we can easily use
substitution to reduce them to standard one-variable quadratic equations, the
Babylonians did not have a way to solve those standard quadratics.  However,
they did consider the problem of finding the dimensions of a rectangular field
given its semiperimeter and area, and had the key substitution used in our
solution method. This is discussed in Gandz's extensive 150-page study of
quadratic equations \cite{gandz-1937}, as well as in Berriman \cite{berriman},
Burton \cite{burton}, Gandz \cite{gandz-1940}, Katz \cite{katz}, and Robson
\cite{robson}.  The ancient Egyptians also had evidence of work with a two-term
quadratic equation, preserved on scraps of a Middle Kingdom papyrus
\cite{edwards-mideast}.

Ancient Chinese mathematicians had solutions to practical problems of
quadratic nature, such as Problem 20 in Chapter 9 of Jiu Zhang Suan Shu
(The Nine Chapters on the Mathematical Art), which was written over several
centuries and completed around 100 A.D. Practical problems of quadratic
nature continued to be considered by other Chinese mathematicians, such as
the 13th-century Yang Hui. See, e.g., the book in Chinese by Zeng
\cite{zeng} or the book in English by Lam \cite{lam}.

The Greeks had several methods of approaching certain types of quadratic
problems, both algebraic and geometric, as surveyed in Eells \cite{eells}.
Heath's translation \cite{heath} of Diophantus \cite{diophantus} from
around 250 A.D. clearly shows the solution of the core problem of finding
two numbers with given sum and product (Book I Problem 27), using the key
parameterization in terms of the average.

Indian mathematicians also had derived a formula for quadratics. Although
Brahmagupta \cite{brahmagupta} did not discover it himself, one root of the
quadratic formula (without derivation) appears in his writings circa 628
A.D. See, e.g., the translation by Colebrooke \cite{colebrooke} or the
commentary by Sharma et al. \cite{sharma}. A derivation due to Sridhara
from around 900 A.D. appears in Puttaswamy \cite{puttaswamy}.

The Persian mathematician al-Khwarizmi published his influential work
\cite{al-khwarizmi} around 825 A.D., where he abstractly considered and
solved the general form of quadratic equations, without starting from
practical applications. His work split into several cases, because he did
not allow numbers to be negative or zero. Consequently, his formulas did
not produce all roots, although they did produce all roots according to the
standards of what a number was at the time.

As mathematics in Western Europe flourished during the Renaissance, successive
formulations and proofs appeared, from Stevin \cite{stevin} to Vi\`ete
\cite{viete} and Descartes \cite{descartes}, ultimately taking on the modern
form that we know today. In the years since then, new proofs have occasionally
appeared, such as two in \emph{The American Mathematical Monthly:} Heaton
\cite{heaton} in 1896 and Cirul \cite{cirul} in 1937.

After an earlier version of this arXiv preprint circulated across the Internet,
references of more recent similar work were identified. The most similar is
Savage \cite{savage}. His approach essentially overlapped in almost all
calculations, but had a pedagogical difference in choice of sign, factoring in
the form $(x+p)(x+q)$ and negating at the end. Perhaps due to a friendly writing
style, that published article has some reversed directions of implication that
are not formally correct. The directional reversals brought in the same extra
assumption as in Footnote \ref{ftCompSqr} when completing the square, creating
another pedagogical difference. That said, those oversights can easily be
corrected by using language similar to our presentation. Gowers \cite{gowers}
also had happened upon a similar approach, while informally presenting a natural
way to deduce the cubic formula.  As he was writing for a different purpose, his
version as written uses Vi\`ete's sum and product relations at the outset,
requiring initial knowledge that there always exist two roots (a pedagogical
difference for first-time Algebra learners), and deducing forward. It can easily
be converted to avoid this existence assumption by using factoring as in our
approach. 

In summary, the author has not yet found a previously-existing book or paper
which states the same pedagogical method as this present article and precisely
justifies the steps, but there exist independent references that contain the
key ideas and can be adapted to achieve this. That said, it is entirely
possible that the method in this present article was previously observed by
people who did not share their findings.

\subsection{Why not centuries ago?}

The two main components of our derivation have existed for hundreds of years
(polynomial factoring converse of Vi\`ete's relations) and for thousands of
years (Babylonian solution to the sum-product problem). Furthermore, the
reduction from the Babylonian problem to a standard quadratic equation has been
well-known for an extremely long time.  Even al-Khwarizmi \cite{al-khwarizmi},
after abstractly analyzing general quadratic equations, showed how to use his
formula to find two numbers with sum 10 and product 21. Like many students in
the modern era, he used substitution to reduce the problem to a single-variable
quadratic equation, and solved it with the quadratic formula. Why, then, didn't
early mathematicians just reverse their steps and find our simple derivation?

Perhaps the reason is because it is actually mathematically nontrivial to
attempt to factor $x^2 + Bx + C = (x - R)(x - S)$ over complex numbers.  Even if
the original quadratic polynomial has real coefficients, it is sometimes
impossible to find two real numbers with sum $-B$ and product $C$.  Early
mathematicians did not know how to reason with a full (algebraically closed)
system of numbers. Indeed, al-Khwarizmi did not even use negative numbers, nor
did Vi\`ete, not to mention the complex numbers that might arise in general.
Perhaps, by the time our mathematical sophistication had advanced to a
sufficient stage, the Babylonian trick had faded out of recent memory, and we
already found the method of completing the square to be sufficiently elementary
for integration into mainstream curriculum.

It is worth noting that the author discovered the solution method in this paper
in the course of filming mathematical explanations, to explain advanced concepts
to particularly young students. Given his audience, he was systematically going
through the middle school math curriculum, creating alternative explanations in
elementary language. To prepare students for the mindset of factoring, he posed
a standalone sum-product problem, designed to be solved via guess-and-check.
While teaching it one evening, his background in coaching math competition
students led him to independently reinvent the Babylonian parameterization in
terms of the average, and to recognize the difference of squares. Later, when
teaching factoring, he suddenly realized that the same technique worked in
general, leading to a simple proof of the quadratic formula! May this story
encourage the reader to think afresh about old things; seeing as how progress
was made on this 4,000 year old topic, more surprises certainly await the light
of discovery.

\end{document}